\documentclass{amsproc}

\usepackage{amsfonts,amsmath,amsxtra,amsthm}
\input xy
\xyoption{all}
\usepackage[all,knot]{xy}

\newtheorem{theorem}{Theorem}[section]
\newtheorem{lemma}[theorem]{Lemma}

\theoremstyle{definition}
\newtheorem{definition}[theorem]{Definition}
\newtheorem{example}[theorem]{Example}
\newtheorem{conjecture}[theorem]{Conjecture}
\newtheorem{corollary}[theorem]{Corollary}

\theoremstyle{remark}
\newtheorem{remark}[theorem]{Remark}

\theoremstyle{proposition}
\newtheorem{proposition}[theorem]{Proposition}

\theoremstyle{propositions}
\newtheorem{propositions}[theorem]{Propositions}
\numberwithin{equation}{section}

\title{$q,t$-Catalan numbers and knot homology}

\author{E. Gorsky}

\address{Department of Mathematics, Stony Brook University, Stony Brook, New York 11794}

\email{egorsky@math.sunysb.edu}

\thanks{Partially supported by the grants RFBR-08-01-00110-a, RFBR-10-01-00678, NSh-8462.2010.1  and the Dynasty fellowship for young scientists.}


\date{}


\begin{document}
\maketitle

\begin{abstract}
We propose an algebraic model of the conjectural triply graded 
homology of S. Gukov, N. Dunfield and J. Rasmussen for some torus knots.
It turns out to be related to the q,t-Catalan numbers of A. Garsia and M. Haiman.
\end{abstract}

\section{Introduction}

In \cite{gaha} A. Garsia and M. Haiman constructed a series of bivariate polynomials $C_{n}(q,t)$.
In \cite{haim} M. Haiman proved that these polynomials have non-negative integer coefficients, and they generalize two known one-parametric deformations of the Catalan numbers, in particular, the value $C_{n}(1,1)$ equals to the $n$-th Catalan number. One of deformations can be expressed in terms of $q$-binomial coefficients, while the second one counts Dyck paths weighted by the area above them. M. Haiman also related these invariants to the geometry of the Hilbert scheme of points on $\mathbb{C}^2$.

Let $Hilb^{n}(\mathbb{C}^2)$ denote the Hilbert scheme of $n$ points on $\mathbb{C}^2$, and let $Hilb^{n}(\mathbb{C}^2,0)$ parametrize 0-dimensional subschemes of length $n$ supported at the origin. Let $V$ be the tautological $n$-dimensional bundle over $Hilb^{n}(\mathbb{C}^2)$.

\begin{theorem}(\cite{haim}, Theorem 2) Consider the natural torus action on $\mathbb{C}^2$ and extend it to Hilbert schemes.
Then $$C_{n}(q_1,q_2)=\chi^{T}(Hilb^{n}(\mathbb{C}^2,0),\Lambda^{n}V),$$
where $q_1$ and $q_2$ are equivariant parameters corresponding to the torus action.
\end{theorem}

We construct a sequence of the bigraded subspaces in the space of symmetric polynomials
such that their Hilbert functions coincide with $C_{n}(t,q)$ for $n\le 4$. 

Let $\Lambda$ denote the ring of symmetric polynomials in the infinite number of variables.
Let $e_k$ denote the elementary symmetric polynomials and $h_k$ denote the complete symmetric polynomials.
One can equip $\Lambda$ with the pair of gradings - one of them is the usual (homogeneous) degree, and the second one is the degree
of a symmetric polynomial as a polynomial in variables $e_k$. In other words, 
$$
S(e_{\alpha_1}\ldots e_{\alpha_r})=\alpha_1+\ldots+\alpha_r,\quad b(e_{\alpha_1}\ldots e_{\alpha_r})=r.
$$
We also define the sequence of spaces $\Lambda(n,r)\subset \Lambda$ which are generated by the monomials with $b$-grading less than or equal to $r$ and $S$-grading equal to $n$.

\begin{definition}
Let $L_n\subset \Lambda$ be the subspace generated by all monomials $h_{\alpha_1}h_{\alpha_2}\ldots h_{\alpha_n}$
such that $\alpha_k\le k$ for all $k$. 
\end{definition}

\begin{theorem}
\label{Ln}
For $n\le 4$, the bivariate Hilbert function of $L_n$ equals to 
\begin{equation}
\sum_{m,r=0}^{\infty}q^{n}t^{r}\dim [(L_n\cap \Lambda(m,r))/(L_n\cap \Lambda(m,r-1))]=q^{n(n-1)/2}C_{n}(q^{-1},t).
\end{equation}
\end{theorem}

The construction of the spaces $L_n$ is expected to be related to some constructions in knot theory.

\begin{definition}(\cite{homfly})
The HOMFLY polynomial $\overline{P}$ is defined by the following skein relation:

\knottips{TF}

$$a\overline{P}\bigl<
\xygraph{
!{0;/r1.0pc/:}
[u(0.5)]
!{\xunderv}
}
\bigr>-a^{-1}\overline{P}\bigl<
\xygraph{
!{0;/r1.0pc/:}
[u(0.5)]
!{\xoverv}
}
\bigr>=(q-q^{-1})\overline{P}
\bigl<
\xygraph{
!{0;/r1.0pc/:}
[u(0.5)]
!{\xunoverv}
}
\bigr>,$$
the multiplication property $\overline{P}(K_1\sqcup K_2)=\overline{P}(K_1)\overline{P}(K_2)$ and 
its non-vanishing at the unknot. One can check that $\overline{P}(unknot)=(a-a^{-1})/(q-q^{-1})$,
and we will also use the {\it reduced} HOMFLY polynomial $P$
$$P(K)(a,q)=\overline{P}(K)(a,q)/\overline{P}(unknot).$$
\end{definition}

The HOMFLY polynomial
unifies the quantum $sl(N)$ polynomial invariants of K $$\overline{P}_N(K)(q)=\overline{P}(K)(a = q^N,q).$$ The original Jones polynomial $J(K)$ equals to $\overline{P}_2(K).$ The HOMFLY
polynomial encodes the Alexander polynomial as well: $\Delta(q)=\overline{P}(K)(a=1,q)$.

The structure of the HOMFLY polynomial for torus knots was described by V. Jones in \cite{jones}. In particular, this result gives the answers for the Alexander, Jones and $sl(N)$ polynomials for all torus knots.
 
More recently, several knot homology theories had been developed: P. Ozsv\'ath and Z. Szab\'o constructed (\cite{os2}) the Heegard-Floer knot homology theory categorifying the Alexander polynomial by the methods of the symplectic topology.  
For all algebraic (and hence torus) knots they managed (\cite{plum}, see also \cite{hedden}) to calculate explicitly the Heegard-Floer homology. It can be reconstructed by a certain combinatorial procedure from the Alexander polynomial.

M. Khovanov (\cite{kho}) constructed a homology theory categorifying the Jones polynomial. Later Khovanov and Rozansky gave a unified construction (\cite{khoro1}) of the homology theories categorifying $sl(N)$  Jones polynomials, and also another homology theory (\cite{khoro2}) categorifying the HOMFLY polynomial.

Although the complexes in the homology theories of Khovanov and Rozansky are defined combinatorially in terms of the knot diagrams, the explicit Poincar\'e polynomials for the corresponding homology groups of torus knots are known only in some particular cases. 

To get all these theories together, Dunfield, Gukov and Rasmussen conjectured (\cite{dgr}) that all these theories are parts, or specializations of a unified picture. Namely, for a given knot $K$ they conjectured the existence of a triply-graded knot homology theory $\mathcal{H}_{i,j,k}(K)$ with the following properties:

\begin{itemize}
{\item {\bf Euler characteristic}. Consider the Poincar\'e polynomial $$\mathcal{P}(K)(a,q,t)=\sum a^{i}q^{j}t^{k}\dim \mathcal{H}_{i,j,k}.$$ Its value at $t=-1$ equals to the value of the reduced HOMFLY polynomial of the knot $K$:
$$\mathcal{P}(K)(a,q,-1)=P(K)(a,q).$$}

{\item {\bf Differentials.} There exist a set of anti-commuting differentials $d_{j}$ for $j\in \mathbb{Z}$ acting in $\mathcal{H}_{*}(K)$. For $N>0$, $d_N$ has triple degree $(-2,2N,-1)$, $d_0$ has degree $(-2,0,-3)$ and for $N<0$ $d_N$ has degree $(-2,2N, -1+2N)$}

{\item {\bf Symmetry.} There exists a natural involution $\phi$ such that
$$\phi d_{N}=d_{-N}\phi$$
for all $N\in \mathbb{Z}$.}

\end{itemize}

For $N\ge 0$, the homology of $d_N$ are supposed to be tightly related to the $sl(N)$ Khovanov-Rozansky homology. Namely, let
$$\mathcal{H}^{N}_{p,k}(K)=\oplus_{iN+j=p}\mathcal{H}_{i,j,k}(K).$$

\begin{conjecture}(\cite{dgr}).
\label{sup}
There exists a homology theory with above properties such that for all $N>1$ the homology of $(\mathcal{H}^{N}_{*}(K),d_N)$ is isomorphic to the $sl(N)$ Khovanov-Rozansky homology. For $N=0$, $(\mathcal{H}^{0}_{*}(K),d_0)$ is isomorphic to the Heegard-Floer knot homology.
The homology of $d_1$ are one-dimensional.
\end{conjecture}

In \cite{rasmussen} J. Rasmussen proved a weaker version of this conjecture. Namely, for all $N>0$ he constructed explicit spectral sequences starting from the Khovanov-Rozansky categorification of HOMFLY polynomial and converging to $sl(N)$ homology.  For the Heegard-Floer homology no relation to the other knot homology theories is known yet.

We propose a conjectural algebraic construction of vector spaces $\mathcal{H}(T_{n,n+1})$ associated with the $(n,n+1)$ torus knots for $n\le 4$. In order to approach the Conjecture \ref{sup}, we prove the following

\begin{theorem}
\label{h}
For $n\le 4$ the Euler characteristic of $\mathcal{H}(n,n+1)$ coincides  with the HOMFLY polynomial of the $(n,n+1)$ torus knot. One can define the differentials $d_0$, $d_1$ and $d_2$
 such that the following properties hold:
\begin{itemize}
\item[1.]{The homology of $\mathcal{H}(n,n+1)$ with respect to the differential $d_1$ is one-dimensional.}
\item[2.]{The homology with respect to $d_2$ is isomorphic to the reduced Khovanov homology of the corresponding knot}
\item[3.]{The homology with respect to $d_0$ is isomorphic to the Heegard-Floer homology of the corresponding knot.}
\end{itemize}
\end{theorem}

Two latter statements are based on the tables from \cite{katlas} and explicit description of the Heegard-Floer homology of algebraic knots proposed in 
\cite{os2} (see also \cite{hedden},\cite{intfun}).

The paper is organized as follows. Section 2 is devoted to the combinatorics of $(q,t)$-Catalan numbers and their polynomial "categorifications". In  Subsection 2.1 we define these numbers and list some of their properties following A. Garsia and M. Haiman.
In  Subsection 2.2 we define the bounce statistic introduced by J. Haglund (\cite{haglund})  and propose a "slicing" construction
dividing a Young diagram into smaller "stable" subdiagrams. In the next subsection we associate a Schur polynomial to a stable Young diagram, and the product of such polynomials for "stable slices" to unstable one. This construction associates a symmetric polynomial to a Dyck path in the $n\times n$ square. It turns out that the subspace generated by these polynomials coincides with the space $L_n$ (defined above), and the gradings of the polynomials are clearly expressed via the area and bounce statistics. This proves  Theorem \ref{Ln}. In  Subsection 2.4 we discuss a generalization of this construction applied to the $(q,t)$-deformation of Schr\"oder numbers defined by J. Haglund.   

Section 3 deals with the HOMFLY polynomials of torus knots and its conjectural categorification. Using the formula of V. Jones, we prove that the coefficients of the HOMFLY polynomial in the power expansion in the variable $a$
can be expressed via certain products of the $q$-binomial coefficients. These coefficients are equal to the generalized Catalan and Schr\"oder numbers. Moreover, the categorification procedure introduces one additional parameter $t$ in the picture, so the resulting coefficients at given powers of $a$ should be some bivariate deformations of the Catalan and Schr\"oder numbers.

Therefore it is quite natural to relate them to the above constructions. Namely, we identify the space $\mathcal{H}(T_{n,m})$ corresponding to a torus knot with a certain subspace in a free polynomial algebra with even and odd generators. This space is equipped with the three gradings: two of them are defined on the ring of symmetric functions as above, and the third one equals to the degree in the skew variables. The differentials of Gukov-Dunfield-Rasmussen are supposed to be
realized as certain differential operators acting on the skew variables. Moreover, we consider the bigger algebra $\mathcal{A}_{n,m}$ acting on $\mathcal{H}(T_{n,m})$. We suppose that for $(n,n+1)$ torus knots the 
space $\mathcal{H}(T_{n,n+1})$ is generated by the volume form and the action of the algebra $\mathcal{A}_{n,n+1}$.
The generators of $\mathcal{A}_{n,n+1}$ can be naturally labelled by the diagonals of the $(n+2)$-gon.

In the Subsection 3.2 we also discuss the "stable limit" of the homology of $(n,m)$-torus knots at $m\rightarrow \infty$, following \cite{dgr}. We identify this limit with the free supercommutative  algebra $\mathcal{H}_n$ with $n-1$ even and $n-1$ odd generators, and show that the grading conditions  define some differentials completely. We compare the resulting constructions and the homology with \cite{katlas} and \cite{dgr}. 

As a byproduct of the above conjectures, we propose an interesting combinatorial conjecture on the limit $q=1$ in triply graded homology. It is well-known in the theory of Heegard-Floer homology that there is a spectral sequence starting from the homology of a given knot and converging to the one-dimensional Heegard-Floer homology of 3-sphere. This means that for any knot the value of the Poincar\'e polynomial for Heegard-Floer homology at $q=1$ equals to 1.

For the triply graded theory, the limit of the Poincar\'e polynomial at $q=1$ is a polynomial in $a$ and $t$.
\begin{conjecture}
Consider the $n\times m$ rectangle and the diagonal in it. Let $D_{n,m}(k)$ denote the set of lattice paths above the diagonal in this rectangle with $k$ marked external corners. For a path $\pi\in D_{n,m}(k)$ let $S(\pi)$
denote the area above $\pi$. Let
$$Q_{n,m}(a,t)=\sum_{k}\sum_{\pi\in D_{n,m}(k)}a^{2k}t^{k+2S(\pi)}.$$
Then the polynomial $Q_{n,m}$ coincides with the limit of the Poincar\'e polynomial for reduced triply graded homology of the torus $(n,m)$-knot at $q=1$.
\end{conjecture}

One can say that the "homological grading" $t$ is related to the area statistics.
This conjecture seems to be coherent to some concepts in mathematical physics (e.g. \cite{gukov}) relating knot homology theories to the geometry of Hilbert schemes and Donaldson-Thomas invariants. 

\bigskip

In \cite{oblomkov} A. Oblomkov and V. Shende developed a remarkable algebro-geometric description of the HOMFLY polynomial for algebraic knots as certain integrals with respect to the Euler characteristic. They proved their conjectures in full generality for all torus knots, and the reduction to the Alexander polynomial for all algebraic knots. We plan to find out the relation between their approach and the one presented here.

The author is grateful to S. Gusein-Zade, S. Gukov, B. Feigin, S. Loktev, A. Gorsky, M. Bershtein and especially to M. Gorsky for lots of useful discussions and remarks. The author also thanks the University of Kumamoto where the part of the work has been done and personally S. Tanabe for the hospitality. The research was partially benefited from the 
support of the "EADS Foundation Chair in mathematics".

The author is grateful to F. Bergeron for pointing that  Theorem \ref{Ln} is not true for $n\ge 5$.

\section{Bivariate Catalan numbers}
\subsection{Properties}

\begin{definition}
A Dyck path in the $n\times n$ square is a lattice path starting at the origin
$(0, 0)$ and ending at $(n, n)$ consisting of North $N (0, 1)$ and East $E (1, 0)$ steps which never
goes below the line $y = x$.

The number $C_n$ of Dyck paths in the $n\times n$ square equals to the $n$-th Catalan number, i. e. $C_n={1\over n+1}{2n\choose n}.$
\end{definition}

In \cite{gaha} A. Garsia and M. Haiman introduced a remarkable two-parametric deformation of the Catalan numbers.

For a cell $x$ in a Young diagram $\mu$ let $l(x),a(x),l'(x),a'(x)$ denote respectively leg, arm, co-leg and co-arm lengths of $x$. Let $$n(\mu)=\sum_{x\in \mu}l(x),\quad n(\mu')=\sum_{x\in \mu}a(x).$$

\begin{definition}(\cite{gaha})
We set
\begin{equation}
\label{defc}
C_{n}(t,q)=\sum_{|\mu|=n}\frac{t^{n(\mu)}q^{n(\mu')}(1-t)(1-q)(\prod_{x\in \mu\setminus (0,0)}(1-t^{l'(x)}q^{a'(x)}))(\sum_{x\in \mu}t^{l'(x)}q^{a'(x)})}{\prod_{x\in \mu}(1-t^{1+l(x)}q^{-a(x)})(1-t^{-l(x)}q^{1+a(x)})}
\end{equation}
\end{definition}

Garsia and Haiman observed that $C_n(t,q)$ is a polynomial with the non-negative integer coefficients (\cite{gahaim}), and $C_n(1,1)$ equals to the Catalan number $c_n$. The geometric meaning of this bivariate deformation of Catalan numbers is described by the following theorem of Haiman.

Let $Hilb^{n}(\mathbb{C}^2)$ be the Hilbert scheme of $n$ points on $\mathbb{C}^2$, and let $Hilb^{n}(\mathbb{C}^2,0)$ parametrize 0-dimensional subschemes of length $n$ supported at the origin. Let $V$ be the tautological $n$-dimensional bundle over $Hilb^{n}(\mathbb{C}^2)$.

\begin{theorem}(\cite{haim}) Consider the diagonal action of the torus $(\mathbb{C}^{*})^2=T$ on $\mathbb{C}^2$ and extend it to Hilbert schemes.
Then $$C_{n}(q_1,q_2)=\chi^{T}(Hilb^{n}(\mathbb{C}^2,0),\Lambda^{n}V),$$
where $q_1$ and $q_2$ are equivariant parameters corresponding to the torus action.
\end{theorem}

As a corollary, $C_n(q_1,q_2)$ is a symmetric function of the parameters $q_1$ and $q_2$.
Two different specializations of $C_n(q_1,q_2)$ are known. 

We will use the standard notation

$$[k]_q=(1-q^k)/(1-q), [k]_q!=[1]_q[2]_q\cdots[k]_q, {n\choose k}_q=[n]_q!/[k]_q![n-k]_q!.$$

\begin{proposition}(\cite{gaha},\cite{haim}) 
1. The values $C_n(q^{-1},q)$ are related to the deformation of Catalan numbers based on $q$-binomial coefficients: 
\begin{equation}
\label{spec}
 q^{n\choose 2}C_n(q^{-1},q)=\frac{1}{[n+1]_q}{2n\choose n}_q.
\end{equation}
2. The values $C_n(1,q)$ coincide with the Carlitz-Riordan (\cite{rior}) $q$-deformation of the Catalan numbers, which are defined by the recursive equation
$$C_n(q)=\sum_{k=0}^{n-1}q^{k}C_k(q)C_{n-1-k}(q), C_0(q)=1.$$
It is also known (e. g. \cite{haglund}) that $C_n(q)=\sum_{\pi}q^{S(\pi)},$ where the summation is done over the set of Dyck paths and $S(\pi)$ denotes the area above the path $\pi$.
\end{proposition}

\subsection{Bounce and area statistics}
 
\begin{definition}
For a Dyck path $\pi$ in the $n\times n$ square we define two statistics, following J. Haglund (\cite{haglund}).
First, $S(\pi)$ is the area above the path $\pi$.

The second one is called the {\it bounce} statistic. Consider a ball starting from the NE corner $(n,n)$.
A ball rolls west until it meets $\pi$, then turns south until it meets the diagonal, then reflects from the diagonal and moves west etc. It finishes at the last point $(0,0)$. During its motion the ball touches the diagonal at points $(j_1,j_1),(j_2,j_2),\ldots$. We define
$$bounce(\pi)=j_1+j_2+\ldots.$$
The ball's path is called {\it bounce path}.

In the below picture the Dyck path is bold and its bounce path is dashed.
\end{definition}

\unitlength=1pt
\begin{picture}(122,300)(50,-50)
\put(101,1){\line(1,0){240}}
\put(341,1){\line(0,1){240}}
\put(101,241){\line(1,0){240}}
\put(101,1){\line(0,1){240}}
\put(101,1){\line(1,1){240}}

{\thicklines
\linethickness{1.5 pt}
\put(221,211){\line(0,1){30}}
\put(161,91){\line(0,1){60}}
\put(161,151){\line(1,0){30}}
\put(191,151){\line(0,1){60}}
\put(191,211){\line(1,0){30}}
\put(101,91){\line(1,0){60}}
\put(101,1){\line(0,1){90}}
\put(221,241){\line(1,0){120}}
}

\put(86,121){\line(1,0){70}}
 
\put(221,204){\line(0,1){5}}
\put(221,197){\line(0,1){5}}
\put(221,190){\line(0,1){5}}
\put(221,183){\line(0,1){5}}
\put(221,176){\line(0,1){5}}
\put(221,169){\line(0,1){5}}
\put(221,162){\line(0,1){5}}
\put(221,155){\line(0,1){5}}
\put(221,148){\line(0,1){5}}
\put(221,141){\line(0,1){5}}
\put(221,134){\line(0,1){5}}
\put(221,127){\line(0,1){5}}
\put(221,121){\line(0,1){4}}

\put(214,121){\line(1,0){5}}
\put(207,121){\line(1,0){5}}
\put(200,121){\line(1,0){5}}
\put(193,121){\line(1,0){5}}
\put(186,121){\line(1,0){5}}
\put(179,121){\line(1,0){5}}
\put(172,121){\line(1,0){5}}
\put(165,121){\line(1,0){5}}
\put(161,121){\line(1,0){3}}

\put(161,86){\line(0,1){5}}
\put(161,79){\line(0,1){5}}
\put(161,72){\line(0,1){5}}
\put(161,65){\line(0,1){5}}
\put(161,61){\line(0,1){2}}

\put(156,61){\line(1,0){5}}
\put(149,61){\line(1,0){5}}
\put(142,61){\line(1,0){5}}
\put(135,61){\line(1,0){5}}
\put(128,61){\line(1,0){5}}
\put(121,61){\line(1,0){5}}
\put(114,61){\line(1,0){5}}
\put(107,61){\line(1,0){5}}
\put(101,61){\line(1,0){4}}

\put(221,114){$j_1$}
\put(161,54){$j_2$}

\put(141,184){$\pi_1$}
\put(131,104){$\pi_2$}

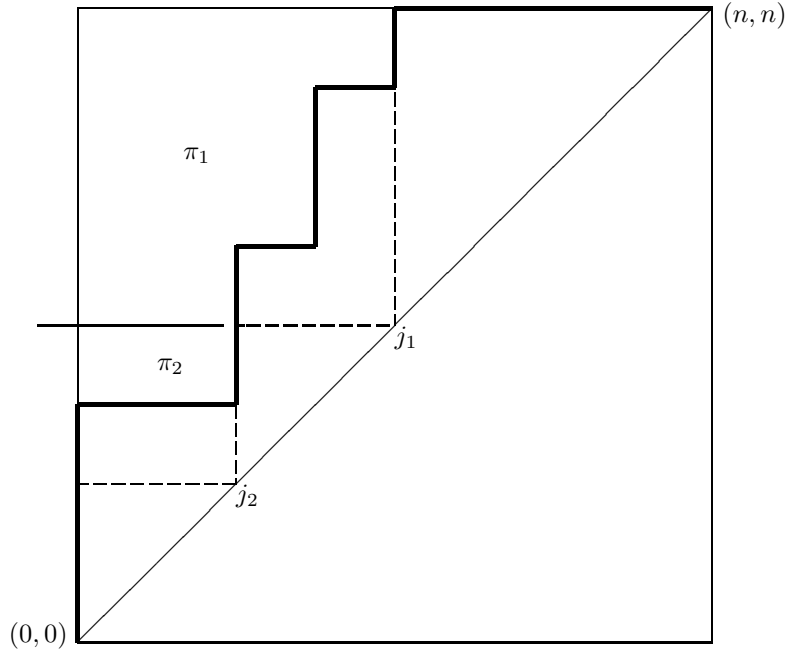
\put(61,-24){Bounce statistic for a Dyck path and slicing of a Young diagram}

\put(75,1){$(0,0)$}
\put(345,236){$(n,n)$}

\end{picture}

It turns out that the area and bounce statistics (which are defined here in a slightly different way than in \cite{haglund}) are related to the $(q,t)$-deformation of the Catalan numbers.

\begin{theorem}(\cite{gahaim},\cite{haglund}) The polynomial $C_n(q_1,q_2)$ can be presented as the following sum over Dyck paths:
\begin{equation}
\label{bou}
C_n(q_1,q_2)=\sum_{\pi} (q_1)^{{n\choose 2}-S(\pi)}(q_2)^{bounce(\pi)}.
\end{equation}
\end{theorem}

\begin{definition}
Consider a Dyck path $\pi$ in the square $n\times n$ and the bounce path for it. Let us continue horisontal bounce lines and cut the Young diagram above $\pi$ along these lines. If bounce points are at $j_1>j_2>\ldots>j_r$, then we get $r$ Young diagrams $\pi_1, \pi_2,\ldots, \pi_r$ corresponding to proper Dyck paths in the squares $n\times n, j_1\times j_1,\ldots, j_r\times j_r$. We will refer to the decomposition $$\pi=\pi_1\sqcup\ldots\sqcup \pi_r$$  
as to the {\bf slicing} of a diagram $T$.
\end{definition}

\begin{definition}
A Dyck path $\pi$ in the square $n\times n$ is {\bf stable}, if $$width(\pi)+height(\pi)<n.$$
\end{definition}

The following properties of slicing and stable paths follows from the definitions.

\begin{propositions}
1. A path $\pi$ is stable if and only if $$bounce(\pi)=width(\pi).$$

2. If $\pi=\pi_1\sqcup\ldots\sqcup \pi_r$ is a slicing of a diagram $\pi$, then slices $\pi_i$ are stable.
Moreover, $width(\pi_m)=j_{m}$, so 
\begin{equation}
\label{slice}
bounce(\pi)=bounce(\pi_1)+\ldots+bounce(\pi_r).
\end{equation}
\end{propositions}

We conclude that bounce and area statistics can be reconstructed from the slicing of the initial path.

\subsection{Symmetric polynomials}

Let $\Lambda$ denote the ring of symmetric polynomials in the infinite number of variables.
Let $e_k$ denote the elementary symmetric polynomials and let $h_k$ denote the complete symmetric polynomials.
One can equip $\Lambda$ with the pair of gradings - one of them is usual degree, and the second one is the degree
of a symmetric polynomial as a polynomial in variables $e_k$. In other words, 
$$
a(e_{\alpha_1}\ldots e_{\alpha_r})=\alpha_1+\ldots+\alpha_r,\quad b(e_{\alpha_1}\ldots e_{\alpha_r})=r.
$$
We also define the sequence of spaces $\Lambda(n,r)\subset \Lambda$ which are generated by the monomials with $b$-grading less than or equal to $r$ and $a$-grading equal to $n$.

We are ready to associate a symmetric polynomial from the ring $\Lambda$ to a Dyck path $\pi$. For stable diagrams the result will not depend of $n$, while for unstable ones it depends on the slicing (and hence on $n$).

\begin{definition}
Let $\pi$ be a stable Young diagram in the square $n\times n$, let $\pi^{*}$ be its transpose, $\pi=(\mu_1,\ldots,\mu_s), \pi^{*}=(\lambda_1, \ldots, \lambda_r)$. We define
the corresponding symmetric polynomial as a Schur polynomial of  $\pi^{*}$.
\begin{equation}
\label{schur}
Z(\pi)=\det(e_{\lambda_{i}-i+j+1})=\det(h_{\mu_{i}-i+j+1}).
\end{equation}
\end{definition}

It is clear that $$a(Z(\pi))=S(\pi), b(Z(\pi))=width(\pi)=bounce(\pi).$$ 

\begin{definition}
Let $\pi$ be a Dyck path in the square $n\times n$, $\pi=\pi_1\sqcup\ldots\sqcup \pi_r$ is its slicing.
Then we define
$$Z(\pi)=Z(\pi_1)\cdot\ldots\cdot Z(\pi_r).$$
\end{definition}

From the equation (\ref{slice}) it follows that the map $Z$ respects both gradings: 
$$a(Z(\pi))=S(\pi), b(Z(\pi))=bounce(\pi).$$ 

The subtle point is that the polynomials $Z(T)$ (as well as Schur polynomials) are homogeneous in the $a$-grading, but not homogeneous in the $b$-grading.
 
{\bf Proof of Theorem \ref{Ln}}.
First, let us remark that for a Dyck path $\pi$ the polynomial $Z(\pi)$ belongs to the space $L_n$.
Consider the slicing $\pi=\pi_1\sqcup\ldots\sqcup \pi_r$.
Using the second equation of (\ref{schur}), one can represent $Z(\pi_j)$ as a determinant in $h'$s whose size is the number of rows in $\pi_j$. By the multiplication of such determinants we get a determinant of a block-diagonal matrix of size $n\times n$. From the definition of the bounce path one can prove that all monomials in the expansion of this determinant belong to the space $L_n$, therefore $Z(\pi)\in L_n$.

Second, all polynomials $Z(\pi)$ are linearly independent since their $h$-lex-minimal terms are different.

Since the dimension of the subspace generated by $Z(\pi)$ equals to the Catalan number, we conclude that $Z(\pi)$ form a basis in $L_n$. 

Moreover, one can check that for $n\le 4$ the $b$-maximal parts of $Z(\pi)$ are linearly independent.  

We conclude that the images of the polynomials $Z(\pi)$ form a basis in all quotients
$[(L_n\cap \Lambda(m,r))/(L_n\cap \Lambda(m,r-1))]$.
Now the statement follows from the equation (\ref{bou}).
$\square$

\begin{remark}
For $n=5$ this proof fails, since the $b$-maximal parts of $Z(\pi)$ are no longer linearly independent
for different $T$. Namely, one can check that $$Z(4,2,2)=h_4(h_2^2-h_1h_3),\quad Z(3,2,2,1)=(h_3h_2-h_1h_4)h_2h_1.$$ 
Both polynomials have gradings $S=8, b=6$, but their $b$-maximal parts are proportional to $e_1^4(e_2^2-e_1e_3)$.
\end{remark}

\subsection{Schr\"oder numbers}

\begin{definition}

A Schr\"oder path is a lattice path starting at the origin
$(0, 0)$ and ending at $(n, n)$ consisting of North $N (0, 1)$, East $E (1, 0)$ and Diagonal $D (1, 1)$ steps which never
goes below the line $y = x$.

\end{definition}

We will denote by $S_{n,k}$ the number of Schr\"oder paths in a square $n\times n$ with exactly $k$ diagonal steps (large Schr\"oder number),
and by $R_{n,k}$ the number of such paths with no $D$ steps on the diagonal $y=x$ (little Schr\"oder number).
It is well-known that $R_{n,k}$ equals to the number of ways to draw $n-k-1$ non-intersecting diagonals in a convex $n$-gon, that is, to the number of $k$-dimensional faces of the
 associahedron.
The combinatorial formula for these numbers looks as (\cite{haglund},\cite{hagl2})

$$S_{n,k}={(2n-k)!\over (n-k+1)!(n-k)!k!}, R_{n,k}={(2n-k)!\over n(n+1)\cdot k!(n-k)!(n-k-1)!}.$$

In \cite{haglund} (see also \cite{can}) a certain bivariate deformation of Schr\"oder numbers was proposed. 
To any Schr\"oder path $\pi$ we associate a Dyck path $T(\pi)$ which is nothing but $\pi$ with all $D$ steps thrown away.

\begin{definition}
Let $S(\pi)$ be the area above the path $\pi$.

Now we define the {\it bounce} statistic. First, consider the Dyck path $T(\pi)$ and the bounce path corresponding to it. 
Let us call the vertical lines of ball's motion {\it peak lines}.
For a $D$-type step $x\in \pi$ let $nump(x)$ denote the number of peak lines to the east from $x$.
Now let
$$b(\pi)=bounce(T(\pi))+\sum_{x}nump(x),$$
where summation is done over all $D$-steps $x$.
\end{definition}

\begin{definition}(\cite{EHKK},\cite{haglund})
The $(q,t)$-Schr\"oder polynomials are defined as
\begin{equation}
S_{n,k}(q,t)=\sum_{\pi}q^{{n\choose 2}+{k\over 2}-S(\pi)}t^{b(\pi)},
\end{equation}
where the summation is done over all $(n,k)$-Schr\"oder paths.
The definition of $R_{n,k}(q,t)$ is analogous.
\end{definition}

It is conjectured (\cite{haglund}) that the polynomials $S_{n,k}(q,t)$ are symmetric in $q$ and $t$.
In what follows we will use the following 

\begin{proposition}(Corollary 4.8.1 in \cite{haglund})
For $0\le k\le n$ 
\begin{equation}
\label{snk}
q^{{n\choose 2}-{k\choose 2}}S_{n,k}(q,q^{-1})={1\over [n-k+1]_{q}}{2n-k\choose n-k, n-k, k}_{q}={[2n-k]!_{q}\over [n-k+1]!_{q}[n-k]!_{q}[k]!_q}
\end{equation}
\end{proposition}

\begin{remark}
The equation (\ref{snk}) can be rewritten as
\begin{equation}
\label{snk2}
\sum_{\pi}q^{S(\pi)+b(\pi)}=q^{k^2\over 2}{1\over [n-k+1]_{q}}{2n-k\choose n-k, n-k, k}_{q}.
\end{equation}
\end{remark}

We conjecture the following analogue of this identity for little Schr\"oder numbers.

\begin{conjecture}
\begin{equation}
\label{rnk}
q^{{n\choose 2}-{k\choose 2}}R_{n,k}(q,q^{-1})={[2n-k]!_{q}\over [n]_q[n+1]_q[n-k-1]!_{q}[n-k]!_{q}[k]!_q}
\end{equation}
\end{conjecture}

Let us introduce the natural analogue for the Schr\"oder paths.

\begin{definition}
Let $\pi$ be a Schr\"oder path, and $T(\pi)$ is a Dyck path defined as above. Consider a slicing of $T(\pi)$, and lift the horisontal cuts of the slicing to the diagram of $\pi$. Now, take away all horisontal lines ending by $D$ steps.
We will receive a Young diagram sliced analogously to $T(\pi)$, which will be called sliced Young diagram of $\pi$.
\end{definition}

We expect the existence some analogue  map $Z$ for the Schr\"oder diagrams. 
Such a map is supposed to take values not in the ring $\Lambda$ itself, but in its extension 
$\Lambda<\xi_1,\xi_2,\ldots>$, where $\xi_j$ are some additional anti-commuting variables. The gradings $S$
and $b$ are extended to these skew variables by the formula
$$S(\xi_j)=j-{1\over 2},b(\xi_j)=0.$$

\section{Homological knot invariants}

\subsection{HOMFLY polynomial for torus knots}

Let $T_{n,m}$ be a torus knot of type $(n,m)$, where $n$ and $m$ are coprime integers, $n<m$.
The explicit expression for the HOMFLY polynomials $P(T_{n,m})$ was found by Jones in \cite{jones}, we'll use it in a slightly rewritten form of \cite{dgr}:
\begin{equation}
P(T_{n,m})=(aq)^{(n-1)(m-1)}{1-q^{-2}\over 1-q^{-2n}}\sum_{b=0}^{n-1}q^{-2mb}(\prod_{i=1}^{b}{a^2q^{2i}-1\over q^{2i}-1})(\prod_{j=1}^{n-1-b}{a^2-q^{2j}\over 1-q^{2j}}).
\end{equation}
To compare the knots with different $n$ and $m$, it is more convenient to get rid of negative powers and consider the rescaled version of this polynomial, namely
\begin{equation}
\label{ps}
P_s(T_{n,m})=(a^{-1}q)^{(n-1)(m-1)}P(T_{n,m})={1-q^2\over 1-q^{2n}}\sum_{b=0}^{n-1}q^{2m(n-1-b)}(\prod_{i=1}^{b}{a^2q^{2i}-1\over q^{2i}-1})(\prod_{j=1}^{n-1-b}{a^2-q^{2j}\over 1-q^{2j}})
\end{equation}
$$={1-q^2\over 1-q^{2n}}\sum_{b=0}^{n-1}q^{2mb}(\prod_{i=1}^{n-1-b}{a^2q^{2i}-1\over q^{2i}-1})(\prod_{j=1}^{b}{a^2-q^{2j}\over 1-q^{2j}}).$$

In \cite{dgr} the expansion of $P_s$ by the powers of $a$ was carefully studied. For example, the following equation holds (as above, we have $n<m$):
\begin{equation}
\label{level}
P_s(T_{n,m})=\sum_{J=0}^{n-1}a^{2J}P_{s}^{J}(T_{n,m}).
\end{equation}

\begin{theorem}
\label{thj}
The following equation for the coefficients $P_{s}^{J}$ holds:
\begin{equation}
\label{HOMFLY}
P_{s}^{k}(T_{n,m})=(-1)^{k}q^{2{k+1\choose 2}}{[m+n-k-1]_{q^2}!\over [n]_{q^2}[m]_{q^2}[k]_{q^2}![m-k-1]_{q^2}![n-k-1]_{q^2}!}.
\end{equation}
\end{theorem}

The proof of this identity can be found in the Appendix.

\begin{remark}
It is not clear from (\ref{ps}), that the right hand side is symmetric in $m$ and $n$ although it should be so.
The coefficients (\ref{HOMFLY}) reveal this symmetry.
\end{remark}

\begin{corollary}
The terms of  top and low degree have the $q$-binomial presentations:
\begin{equation}
\label{toplow}
P_{s}^{(n-1)}(T_{n,m})={(-1)^{n-1}q^{n(n-1)}\over [n]_{q^2}}{m-1\choose n-1}_{q^2},\\
P_{s}^{0}(T_{n,m})={1\over [n]_{q^2}}{m+n-1\choose n-1}_{q^2}
\end{equation}

At the limit $q=1$ we get

$$P_{s}^{(n-1)}(T_{n,m})(q=1)={(-1)^{n-1}\over n}{m-1\choose n-1},\\
P_{s}^{0}(T_{n,m})(q=1)={1\over n}{m+n-1\choose n-1}.$$

Also an interesting "blow-up" equation follows from (\ref{toplow}):

\begin{equation}
\label{blow} 
P_{s}^{0}(T_{n,m})=(-1)^{n-1}q^{-n(n-1)}P_{s}^{(n-1)}(T_{n,m+n})=(-1)^{m-1}q^{-m(m-1)}P_{s}^{(m-1)}(T_{m,m+n}).
\end{equation}

\end{corollary}

\begin{corollary}
If we focus on the case $m=n+1$, we have
\begin{equation}
\label{n+1}
P_{s}^{k}(T_{n,n+1})=(-1)^{k}q^{k(k+1)}{1\over [n-k]_{q^2}}{n-1\choose k}_{q^2}{2n-k\choose n+1}_{q^2}.
\end{equation}
At the limit {q=1} we have
$$(-1)^{k-1}P_{s}^{k}(T_{n,n+1})={1\over n-k}{n-1\choose k}{2n-k\choose n+1},$$
what is 
equal to the little Schr\"oder number $R_{n,k}$. At the lowest level $k=0$ we get the $n$-th Catalan number. 
\end{corollary}

\begin{definition}
We call  a Dyck path marked if some of its external corners are marked. 
\end{definition}

\begin{theorem}
The number of marked Dyck paths in the rectangle $m\times n$ with $k$ marks equals to
$${(m+n-k-1)!\over m\cdot n\cdot (n-k-1)!(m-k-1)!k!}.$$
\end{theorem}

\begin{proof}
Follows from the Lemmas \ref{Nara1} and \ref{Nara2} from the Appendix.
\end{proof}

\begin{corollary}
\label{dyck}
The coefficient $P_{s}^{k}(T_{n,m})$ of the HOMFLY polynomial for $(n,m)$ torus knot is a certain $q$-deformation
of the number of marked Dyck paths in the rectangle $n\times m$ with $k$ marks.
\end{corollary}

The Corollary \ref{dyck} means that the coefficients at $a^{2k}$ in the Poincar\'e polynomial of the Gukov-Dunfield-Rasmussen homology  of the $(n,m)$-torus knot should be certain $(q,t)$-deformations of the above combinatorial data. For example, we know that $P_{s}^{2k}(T_{n,n+1})$ is a $q$-deformation of the Schr\"oder number $R_{n,k}$, and it is natural to assume that $\mathcal{P}_{s}^{2k}(T_{n,n+1})$ is related to the $(q,t)$-deformation of this number. By (\ref{rnk}) we have
$$q^{2{n\choose 2}-2{k\choose 2}}R_{n,k}(q^2,q^{-2})={[2n-k]!_{q^2}\over [n]_{q^2}[n+1]_{q^2}[n-k-1]!_{q^2}[n-k]!_{q^2}[k]!_{q^2}},$$
$$q^{2{n\choose 2}}R_{n,k}(q^2,q^{-2})=(-1)^{k}q^{-2k}P_{s}^{2k}(T_{n,n+1})=q^{-2k}\mathcal{P}_{s}^{2k}(T_{n,n+1})(q,-1),$$
what motivates the following 

\begin{conjecture}
The following equation holds:
\begin{equation}
\label{pp}
\mathcal{P}_{s}^{2k}(T_{n,n+1})(q,t)=q^{2{n\choose 2}+2k}t^{2{n\choose 2}+3k}R_{n,k}(q^{-2}t^{-2},q^2)=q^{k}t^{2k}\sum_{\pi}q^{2(S(\pi)+b(\pi))}t^{2S(\pi)},
\end{equation}
where summation in the right hand side is done over all $(n,k)$-Schr\"oder paths with no $D$ steps on the diagonal.
\end{conjecture}

\begin{corollary}
The following equation holds:
\begin{equation}
\label{pp}
\mathcal{P}_{s}^{0}(T_{n,n+1})(q,t)=q^{2{n\choose 2}}t^{2{n\choose 2}}C_{n}(q^{-2}t^{-2},q^2)=\sum_{\pi}q^{2(S(\pi)+b(\pi))}t^{2S(\pi)},
\end{equation}
where summation in the right hand side is done over all Dyck paths in $n\times n$ square.
\end{corollary}

Since $(q,t)$-Schr\"oder number are supposed to be symmetric in $q$ and $t$, one can check the symmetry property for $\mathcal{P}$ 
which agrees with the properties of the involution $\phi$ from Conjecture \ref{sup}.

The following equation is a corollary of (\ref{pp}):
$$\mathcal{P}_{s}^{2k}(T_{n,n+1})(1,t)=t^{2k}\sum_{\pi}t^{2S(\pi)},$$
where summation is over all $(n,k)$-Schr\"oder paths with no $D$ steps on the diagonal and $S(\pi)$ denotes the area above the path.
We generalize this remark to the following

\begin{conjecture}
The following equation holds:
\begin{equation}
\label{D1}
\mathcal{P}_{s}^{2k}(T_{n,m})(1,t)=t^{k}\sum_{\pi}t^{2S(\pi)},
\end{equation}
where summation is over all marked Dyck paths in the $m\times n$ rectangle with $k$ marks.
\end{conjecture}

\subsection{Stable limit}

The right hand side of (\ref{ps}) in the limit $m\to \infty$ tends to
$$P_{s}(T_n)=\lim_{m\to\infty}P_{s}(T_{n,m})=\prod_{k=1}^{n-1}{(1-a^2q^{2k})\over (1-q^{2k+2})}.$$

It is natural to consider the behaviour of the Gukov-Dunfield-Rasmussen homology in this limit too.
Following the discussions in Section 6 of \cite{dgr}, we conjecture that the limit homology $\mathcal{H}(T_n)=\lim_{m\to\infty}\mathcal{H}(T_{n,m})$ is a free polynomial algebra
with $n-1$ even generators with gradings $(0,2k+2,2k)$ and $n-1$ odd generators with gradings $(2,2k,2k+1)$,
and therefore 
$$\mathcal{P}_{s}(T_n)=\prod_{k=1}^{n-1}{(1+a^2q^{2k}t^{2k+1})\over (1-q^{2k+2}t^{2k})}.$$

We denote the odd generators by $\xi_1,\ldots,\xi_{n-1}$, and even generators by $e_1,\ldots,e_{n-1}$.
The notation for even generators is motivated by the above constructions related with $(q,t)$-Catalan numbers. To be more precise, we identify $e_k$ with the $k$-th elementary symmetric polynomial, and the even part of $\mathcal{H}(T_{n})$ with the ring of symmetric polynomials in $n-1$ variables. Recall that we had two natural gradings on this ring defined by the equations
$$S(e_k)=k, b(e_k)=1.$$
Therefore the triple grading on the even part equals to $(0,2(S+b),2b)$.

Let us construct the action of the differentials on $\mathcal{H}(T_n)$.
The differentials send $\xi_k$ to some polynomials in $e_m$, and they are extended to the whole algebra by the Leibnitz rule. Taking into account the gradings, one can uniquely guess the equations
$$d_{-n}(\xi_{k})=\delta_{k,n}, d_0(\xi_k)=e_{k-1}, d_1(\xi_k)=e_k.$$
 
Let us compute the homology of $\mathcal{H}(T_n)$ with respect to differentials $d_N$.
From the properties of the Koszul complex one can deduce the following

\begin{propositions}
1. The complexes $(\mathcal{H}(T_n),d_{-N})$ are acyclic.

2. The homology of $(\mathcal{H}(T_n),d_0)$ is the polynomial algebra generated by $\xi_1$ and $e_{n-1}$.

3. The homology of $(\mathcal{H}(T_n),d_1)$ is one-dimensional and generated by 1.
\end{propositions}

The construction of the higher differentials is less restricted by the grading, however for small degrees
one has no choice but to define 
$$d_2(\xi_2)=e_1^2, d_2(\xi_3)=e_1e_2, d_3(\xi_3)=e_1^3.$$

\begin{example}
The homology of $(\mathcal{H}(T_3), d_2)$ is generated by $\xi_1, e_2$ and $e_1$ modulo relation $e_1^2=0$
since $d_2(\xi_2)=e_1^2$. 
These generators have  gradings $(2,2,3),(0,6,4)$ and $(0,4,2)$,
so the Poincar\'e polynomial for these homology equals to
$${(1+q^4t^2)(1+a^2q^2t^3)\over (1-q^6t^4)}$$.
\end{example}

\begin{example}
The homology of $(\mathcal{H}(T_4), d_2)$ is generated by the elements $\xi_1, e_2, e_3, e_1$ and $\mu=e_1\xi_3-e_2\xi_2$ modulo relation $e_1^2=0, e_1e_2=0, e_1\mu=0$, so it is isomorphic to
$$H(\mathcal{H}(T_4),d_2)=\mathbb{C}[\xi_1,e_3]\otimes(<e_1>\oplus \mathbb{C}[\mu,e_2]).$$
The Poincar\'e polynomial for this homology equals to
$${(1+a^2q^2t^3)\over(1-q^8t^6)}[q^4t^2+{1+a^2q^{10}t^9\over 1-q^6t^4}]$$.
\end{example}

One can compare these answers with \cite{dgr}.

\subsection{$(2,k)$ and $(3,k)$ torus knots}

For the torus knots with the small number of strands (2 or 3) it is possible to give a clear algebraic description
of the structure predicted in \cite{dgr}: the triply graded homology and surviving differentials $d_{0}, d_{\pm 1}, d_{\pm 2}$.

\begin{conjecture}
\label{conj23}
The triply graded homology for the $(n,k)$ torus knot for $n\le 3$ can be realized as a subspace in $\mathcal{H}(T_n)$.
This subspace generated by its top level and the action of the differentials. The top level subspace is generated by the following monomials:

$$e_1^{i}\xi_1,\quad 2i\le k-3\qquad \mbox{\rm for}\quad n=2,$$ 
$$e_1^{i}e_2^{j}\xi_1\xi_2,\quad i+3j\le k-4\qquad \mbox{\rm for}\quad n=3.$$ 
\end{conjecture}

\begin{example}
The homology of the trefoil knot $T_{2,3}$ is generated over $d_{\pm 1}$ by one element $\xi_1$. 
The differentials act as $d_{-1}(\xi_1)=1$ and $d_1(\xi_1)=e_1$.
\end{example}

\begin{example} The homology of the trefoil knot $T_{3,4}$ is generated over the differentials $d_{0}, d_{\pm 1}, d_{\pm 2}$ by one element $\xi_1\xi_2$. 
This homology is presented in \cite{dgr} by a diagram with three levels. On the level 2 we have one element
$\xi_1\xi_2$. On the level 1 we have
$$d_{-2}(\xi_1\xi_2)=\xi_1,d_{-1}(\xi_1\xi_2)=\xi_2,d_{0}(\xi_1\xi_2)=e_1\xi_1,d_{1}(\xi_1\xi_2)=e_1\xi_2-e_2\xi_1,d_{2}(\xi_1\xi_2)=e_1^2\xi_1,$$
and on the level 0 we have
$$d_{-1}(\xi_1)=d_{-2}(\xi_2)=1,$$
$$d_1(\xi_1)=d_0(\xi_2)=d_{-1}(e_1\xi_1)=d_{-2}(e_1\xi_2-e_2\xi_1)=e_1;$$
$$d_1(\xi_2)=-d_{-1}(e_1\xi_2-e_2\xi_1)=e_2$$
$$d_2(\xi_2)=d_1(e_1\xi_1)=d_0(e_1\xi_2-e_2\xi_1)=d_{-1}(e_1^2\xi_1)=e_1^2,$$
$$d_2(e_1\xi_2-e_2\xi_1)=d_1(e_1^2\xi_1)=e_1^3.$$
These equations represent exactly the dot diagram for the homology of $T_{3,4}$ presented in \cite{dgr}:

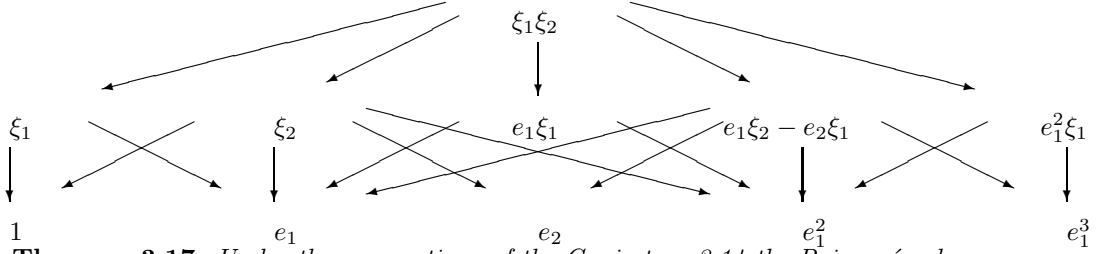
\unitlength=1pt
\begin{picture}(122,100)(50,120)
\put(221,191){$\xi_1\xi_2$}
\put(31,151){$\xi_1$}
\put(131,151){$\xi_2$}
\put(221,151){$e_1\xi_1$}
\put(301,151){$e_1\xi_2-e_2\xi_1$}
\put(421,151){$e_1^2\xi_1$}
\put(31,111){$1$}
\put(131,111){$e_1$}
\put(231,111){$e_2$}
\put(331,111){$e_1^2$}
\put(431,111){$e_1^3$}
\put(201,196){\vector(-2,-1){50}}
\put(196,201){\vector(-4,-1){130}}
\put(261,196){\vector(2,-1){50}}
\put(266,201){\vector(4,-1){130}}
\put(231,186){\vector(0,-1){20}}
\put(31,146){\vector(0,-1){20}}
\put(61,156){\vector(2,-1){50}}

\put(101,156){\vector(-2,-1){50}}
\put(161,156){\vector(2,-1){50}}
\put(166,161){\vector(4,-1){130}}
\put(131,146){\vector(0,-1){20}}

\put(201,156){\vector(-2,-1){50}}
\put(261,156){\vector(2,-1){50}}

\put(301,156){\vector(-2,-1){50}}
\put(296,161){\vector(-4,-1){130}}
\put(361,156){\vector(2,-1){50}}
\put(331,146){\vector(0,-1){20}}

\put(401,156){\vector(-2,-1){50}}
\put(431,146){\vector(0,-1){20}}

\end{picture}

\end{example}

\begin{theorem}
\label{t23}
Under the assumptions of the Conjecture \ref{conj23} the Poincar\'e polynomial for the triply-graded homology
has a form:
\begin{equation}
\label{2k}
\mathcal{P}_{2,2k+1}=\sum_{i=0}^{k}
q^{4i}t^{2i}+a^2\sum_{i=0}^{k-1}q^{4i+2}t^{2i+3},
\end{equation}

\begin{equation}
\label{3k}
\mathcal{P}_{3,k}=1+\sum_{0<i+3j\le k-1}q^{4i+6j}t^{2i+4j}(1+a^2q^{-2}t)+a^{4}q^{6}t^{8}\sum_{i+3j\le k-4}q^{4i+6j}t^{2i+4j}(1+a^{-2}q^{2}t^{-1}).
\end{equation}
 
\end{theorem}

\begin{remark}
The polynomials (\ref{2k}) and (\ref{3k}) coincide with the ones conjectured in \cite{dgr}.
\end{remark}

\begin{proof}
As a vector space, the homology of $T_{2,2k+1}$ is generated by $$1, e_1, e_1^2\ldots e_1^{k}\qquad \mbox{\rm on lower level},$$
$$\xi_1, e_1\xi_1, e_1^2\xi_1\ldots e_1^{k-1}\xi_1\qquad \mbox{\rm on top level}.$$
This implies the equation (\ref{2k}).

By Conjecture \ref{conj23}, we know the structure of $\mathcal{H}(3,k)$ as a vector space at the level 2.
Let us describe it on levels 1 and 0. Remark that $\mathcal{H}(3,k)$ is generated by its top level,
and the differentials commute with the multiplication by even variables, so the levels in $\mathcal{H}(3,k)$ are
just the products of the top level even monomials with the corresponding levels of $\mathcal{H}(3,4)$.

Therefore on the lower level we get 
$$<1,e_1,e_1^2,e_1^3,e_2>\otimes <e_1^{i}e_2^{j}|i+3j\le k-4>= <e_1^{i}e_2^{j}|3i+j\le k-1>.$$
It rests to know that the homology of $d_1$ is one-dimensional and generated by 1, so
we can compute the homology at level 1 as well.
 
\end{proof}

\begin{corollary}
The reduced $sl(2)$ homology of $(2,n)$ torus knot coincides with the triply-graded homology as $d_2=0$.
\end{corollary}

\begin{corollary}
The reduced $sl(2)$ homology of $(3,n)$ torus knot is spanned by the monomials 
\begin{equation}
\label{monomials}
1, e_2^{i}, e_1e_2^{j}, e_2^{j}\xi_1, e_1e_2^{i-1}\xi_1
\end{equation}
with $0<3i\le n-1, 3j\le n-2$.
\end{corollary}

\begin{proof}
Since $d_2(\xi_1)=0, d_2(\xi_2)=e_1^2$, to compute the homology of $\mathcal{H}(T_{3,n})$ with respect to $d_2$ one should 
throw away from $\mathcal{H}(T_{3,n})$ all monomials divisible by  $\xi_2$ and $e_1^2$. Therefore it is spanned by the monomials of the form (\ref{monomials}) belonging to $\mathcal{H}(T_{3,n})$.
\end{proof}

One can compare these results with \cite{katlas}.

\begin{example}
The reduced $sl(2)$ homology of $(3,4)$ torus knot is spanned by the monomials 
$1, e_1, e_2, \xi_1, e_1\xi_1$, what gives the following Poincare polynomial:
$$\mathcal{P}_{2}(T_{3,4})=1+q^4t^2+q^6t^4+q^6t^3+q^{10}t^{5}.$$
\end{example}

\subsection{$(4,5)$ knot}

Following the above description of the triply graded homology for $(2,n)$ and $(3,n)$ torus knots, it is natural to assume that the 
homology of a torus knot is generated by its top level under the action of some algebra. 

\begin{conjecture}
There exists a sequence of algebras $\mathcal{A}_{n,m}$ together with their representations in the stable homology space $\mathcal{H}(T_{n,m})$ satisfying the following conditions:

1) The action of $\mathcal{A}_{n,m}$ commutes with the multiplication by the polynomials in $e_1,\ldots e_n$

2) For every $k$ ($|k|<n$) the differential $d_k$ belongs to $\mathcal{A}_{n,m}$

3) The space $\mathcal{H}(T_{n,m})$ is generated by its top level under the action of $\mathcal{A}_{n,m}$.

4) The top level of homology for $(m,m+n)$ torus knot coincides with the zero level of the homology 
of $(m,n)$ torus knot multiplied by the "volume form" $\xi_1\xi_2\ldots \xi_{m-1}$  (compare with (\ref{blow})).
\end{conjecture}

For $(2,n)$ and $(3,n)$ knots the corresponding algebras are supposed to be generated by the differentials $d_{0}, d_{\pm 1}, d_{\pm 2}$.
Unfortunately, for bigger knots this is not true: for example, $\mathcal{H}_{n,n+1}$ has dimension 1 on the top level, so it is supposed to be generated by the "volume form" $\xi_1\xi_2\ldots \xi_{n-1}$. On the level $(n-2)$ it has dimension $R(n,n-2)=(n+2)(n-1)/2$ that is greater than the number $(2n-1)$ of available differentials. 

For example, for $n=4$ we have $R(4,2)=9$ and 7 differentials at our disposal. Nevertheless, by grading reasons one can uniquely guess the algebraic description of the differentials as well as two missing operators.  Presuming the action of differentials on $\xi_1$ and $\xi_2$ being the same as above, one can extend it to $\xi_3$ by the formula
$$d_{-1}(\xi_3)=d_{-2}(\xi_3)=0,\quad d_{-3}(\xi_3)=1,\quad$$ $$ d_{0}(\xi_3)=e_2,\quad d_{1}(\xi_3)=e_3,\quad d_2(\xi_3)=e_1e_2,\quad d_3(\xi_3)=e_1^3.$$
The two missing operators are 
$$\alpha_1=e_1^{-1}d_3, \alpha_2=e_1^{-2}d_3.$$
Suppose that the space $\mathcal{H}(T_{4,5})$ is generated by the "volume form" $\xi_1\xi_2\xi_3$ under the action of the differentials and the operators $\alpha_1$ and $\alpha_2$.  
Let us describe the basis on each level in this space and indicate the gradings of the basis elements.
For each basis element we also indicate the element in $\mathcal{A}_{4,5}$ producing it from $\xi_1\xi_2\xi_3$.

\begin{tabular}{|c|c|c|c|}
\hline
Level & $\mbox{\rm Basis element}$ & $\mbox{\rm Element of}\quad \mathcal{A}_{4,5}$ & $(a,q,t)$\\
\hline
$3$ & $\xi_1\xi_2\xi_3$ & $1$ & $(6,12,15)$\\
\hline
$2$ & $\xi_1\xi_2$ & $d_{-3}$ & $(4,6,8)$\\
\hline
$2$ & $\xi_1\xi_3$ & $d_{-2}$ & $(4,8,10)$\\
\hline
$2$ & $\xi_2\xi_3$ & $d_{-1}$ & $(4,10,12)$\\
\hline
$2$ & $e_1\xi_1\xi_2$ & $\alpha_2$ & $(4,10,10)$\\
\hline
$2$ & $e_1\xi_1\xi_3-e_2\xi_1\xi_2$ & $d_{0}$ & $(4,12,12)$\\
\hline
$2$ & $e_1\xi_2\xi_3-e_2\xi_1\xi_3+e_3\xi_1\xi_2$ & $d_{1}$ & $(4,14,14)$\\
\hline
$2$ & $e_1^2\xi_1\xi_2$ & $\alpha_1$ & $(4,14,12)$\\
\hline
$2$ & $e_1^2\xi_1\xi_3-e_1e_2\xi_1\xi_2$ & $d_{2}$ & $(4,16,14)$\\
\hline
$2$ & $e_1^3\xi_1\xi_2$ & $d_{3}$ & $(4,18,14)$\\
\hline
$1$ & $\xi_1$ & $d_{-2}d_{-3}$ & $(2,2,3)$\\
\hline
$1$ & $\xi_2$ & $d_{-1}d_{-3}$ & $(2,4,5)$\\
\hline
$1$ & $\xi_3$ & $d_{-1}d_{-2}$ & $(2,6,7)$\\
\hline
$1$ & $e_1\xi_1$ & $d_{0}d_{-3}$ & $(2,6,5)$\\
\hline
$1$ & $e_1\xi_2-e_2\xi_1$ & $d_{1}d_{-3}$ & $(2,8,7)$\\
\hline
$1$ & $e_1\xi_3-e_3\xi_1$ & $d_{1}d_{-2}$ & $(2,10,9)$\\
\hline
$1$ & $e_1\xi_3-e_2\xi_2$ & $d_{0}d_{-1}$ & $(2,10,9)$\\
\hline
$1$ & $e_1^2\xi_1$ & $d_{2}d_{-3}$ & $(2,10,7)$\\
\hline
$1$ & $e_1^3\xi_1$ & $d_{2}\alpha_2$ & $(2,14,9)$\\
\hline
$1$ & $e_1^4\xi_1$ & $d_{2}\alpha_1$ & $(2,18,11)$\\
\hline
$1$ & $e_1^5\xi_1$ & $d_{2}d_{3}$ & $(2,22,13)$\\
\hline
$1$ & $e_1\xi_2$ & $d_{-1}\alpha_2$ & $(2,8,7)$\\
\hline
$1$ & $e_1^2\xi_2$ & $d_{-1}\alpha_1$ & $(2,12,9)$\\
\hline
$1$ & $e_1^3\xi_2$ & $d_{-1}d_{3}$ & $(2,16,11)$\\
\hline
$1$ & $e_1^2\xi_2-e_1e_2\xi_1$ & $d_{1}\alpha_2$ & $(2,12,9)$\\
\hline
$1$ & $e_1^3\xi_2-e_1^2e_2\xi_1$ & $d_{1}\alpha_1$ & $(2,16,11)$\\
\hline
\end{tabular}

\begin{tabular}{|c|c|c|c|}
\hline
Level & $\mbox{\rm Basis element}$ & $\mbox{\rm Element of}\quad \mathcal{A}_{4,5}$ & $(a,q,t)$\\
\hline
$1$ & $e_1^2\xi_3-e_1e_2\xi_2$ & $d_{2}d_{-1}$ & $(2,14,11)$\\
\hline
$1$ & $e_2\xi_3-e_3\xi_2$ & $d_{1}d_{-1}$ & $(2,12,11)$\\
\hline
$1$ & $e_1^4\xi_2-e_1^3e_2\xi_1$ & $d_{1}d_{3}$ & $(2,20,13)$\\ 
\hline
$1$ & $(e_2^2-e_1e_3)\xi_1-e_1e_2\xi_2+e_1^2\xi_3$ & $d_{1}d_{0}$ & $(2,14,11)$\\
\hline
$1$ & $(e_1e_2^2-e_1^2e_3)\xi_1-e_1^2e_2\xi_2+e_1^3\xi_3$ & $d_{1}d_{2}$ & $(2,18,13)$\\
\hline
$0$ & $1$ & $d_{-1}d_{-2}d_{-3}$ & $(0,0,0)$\\
\hline
$0$ & $e_1$ & $d_{1}d_{-2}d_{-3}$ & $(0,4,2)$\\
\hline
$0$ & $e_2$ & $d_{1}d_{-1}d_{-3}$ & $(0,6,4)$\\
\hline
$0$ & $e_3$ & $d_{1}d_{-1}d_{-2}$ & $(0,8,6)$\\
\hline
$0$ & $e_1^2$ & $d_{1}d_{0}d_{-3}$ & $(0,8,4)$\\
\hline
$0$ & $e_1^3$ & $d_{1}d_{2}d_{-3}$ & $(0,12,6)$\\
\hline
$0$ & $e_1^4$ & $d_{1}d_{2}\alpha_2$ & $(0,16,8)$\\
\hline
$0$ & $e_1^5$ & $d_{1}d_{2}\alpha_1$ & $(0,20,10)$\\
\hline
$0$ & $e_1^6$ & $d_{1}d_{2}d_{3}$ & $(0,24,12)$\\
\hline
$0$ & $e_1e_2$ & $d_{1}d_{-1}\alpha_2$ & $(0,10,6)$\\
\hline
$0$ & $e_1e_3-e_2^2$ & $d_{1}d_{0}d_{-1}$ & $(0,12,8)$\\
\hline
$0$ & $e_1^2e_3-e_1e_2^2$ & $d_{1}d_{2}d_{-1}$ & $(0,16,10)$\\
\hline
$0$ & $e_1^2e_2$ & $d_{1}d_{-1}\alpha_1$ & $(0,14,8)$\\
\hline
$0$ & $e_1^3e_2$ & $d_{1}d_{-1}d_{3}$ & $(0,18,10)$\\
\hline
\end{tabular} 

\begin{remark}
The homology of $d_1$ is one-dimensional and spanned by 1. 
\end{remark}

\begin{remark}
The homology of $d_2$ is spanned by 
$$1, e_1, \xi_1, e_2, e_1\xi_1, e_3,  e_2\xi_1, e_1e_3-e_2^2, e_1\xi_3-e_2\xi_2,$$  
what agrees with the generating function for the
reduced Khovanov homology of $T_{4,5}$ (compare with \cite{katlas}):
$$\mathcal{P}_{2}(T_{4,5})=1+q^{4}t^2+q^{6}t^3+q^{6}t^4+q^{10}t^5+q^{8}t^6+q^{12}t^7+q^{12}t^8+q^{14}t^9.$$
\end{remark} 

\begin{remark}
At the lower level of the triply-graded homology, we get the space $L_4$.
In particular, this space has no monomial baisis.
\end{remark}
 
\section{Appendix}

The following lemma is a well known $q$-analogue of the binomial identity.

\begin{lemma}
\begin{equation}
\label{qbinom}
(1+z)(1+qz)\cdot\ldots\cdot(1+q^{n-1}z)=\sum_{j=0}^{n}{n\choose j}_{q}q^{j\choose 2}z^{j}.
\end{equation}
\end{lemma}

\begin{proof}
Induction by $n$.
\end{proof}

{\bf Proof of the Theorem \ref{thj}}
First, by (\ref{qbinom}), we have
$$\prod_{i=1}^{n-1-b}(1-a^2q^{2i})=\sum_{i=0}^{n-1-b}(-1)^{i}a^{2i}q^{2{i+1\choose 2}}{n-1-b\choose i}_{q^{2}},$$
$$\prod_{j=1}^{b}(a^2-q^{2j})=\sum_{j=0}^{b}(-1)^{b-j}a^{2j}q^{2{b-j+1\choose 2}}{b\choose j}_{q^2}.$$
If we multiply these expressions and take the coefficient at $a^{2k}$, we get
$$\sum_{i=k-b}^{n-1-b}(-1)^{b-k}q^{2{i+1\choose 2}+2{b-k+i+1\choose 2}}{n-1-b\choose i}_{q^2}{b\choose k-i}_{q^2}.$$
Therefore 
$$P_{s}^{2k}(T_{n,m})={1\over [n]_{q^2}}\sum_{b=0}^{n-1}q^{2mb}\sum_{i=k-b}^{n-1-b}(-1)^{b-k}q^{2{i+1\choose 2}+2{b-k+i+1\choose 2}}\times$$
$${[b]_{q^2}![n-1-b]_{q^2}!\over [i]_{q^2}![n-1-b-i]_{q^2}![k-i]_{q^2}![b-k+i]_{q^2}!}{(1-q^2)^{n-1}\over [b]_{q^2}![n-1-b]_{q^2}!}=$$
$${(1-q^2)^{n-1}\over [n]_{q^2}}\sum_{i=0}^{k}{(-1)^{i}q^{2m(k-i)}\cdot q^{2{i+1\choose 2}}[k]_{q^2}!\over [n-1-k]_{q^2}![k]_{q^2}![i]_{q^2}![k-i]_{q^2}!}\times $$ $$\sum_{b=k-i}^{n-1-i}(-1)^{b-k+i}q^{2m(b-k+i)}q^{2{b-k+i+1\choose 2}}{[n-1-k]_{q^2}!\over [n-1-b-i]_{q^2}![b-k+i]_{q^2}!}.$$
Now by (\ref{qbinom}) we simplify the inner sum, denoting $l=b-k+i$:
$$\sum_{l=0}^{n-1-k}(-1)^{l}q^{2ml}q^{2{l+1\choose 2}}{n-1-k\choose l}_{q^2}=(1-q^{2m+2})(1-q^{2m+4})\cdot\ldots\cdot(1-q^{2m+2n-2-2k})=$$ $$(1-q^2)^{n-1-k}{[m+n-1-k]_{q^2}!\over [m]_{q^2}!},$$
And this sum does not depend on $i$.
Analogously we have 
$$\sum_{i=0}^{k}(-1)^{i}q^{2{i+1\choose 2}}q^{2m(k-i)}{k\choose i}_{q^2}=(q^{2m}-q^2)\cdot\ldots\cdot(q^{2m}-q^{2k})=$$ $$(-1)^{k}(1-q^2)^{k}q^{2{k+1\choose 2}}{[m-1]_{q^2}!\over [m-k-1]_{q^2}!}.$$
Finally,
$$P_{s}^{2k}(T_{n,m})=(-1)^{k}q^{2{k+1\choose 2}}{[m-1]_{q^2}![m+n-1-k]_{q^2}!\over [n]_{q^2}[m]_{q^2}![m-k-1]_{q^2}![n-1-k]_{q^2}![k]_{q^2}!}.$$
$\square$

To what follows we will need the the following sequence of "generalized Narayana numbers" (for given $n,m$):
$$N_{k}={(m-1)!(n-1)!\over k!(k+1)!(m-k-1)!(n-k-1)!}={1\over k+1}{m-1\choose k}{n-1\choose k}.$$

\begin{lemma}
\label{Nara1}
$N_{k}$ equals to the number of Dyck paths in $m\times n$ rectangle with $k$ external corners.
\end{lemma}

\begin{proof}
First, let us remark that ${m-1\choose k}{n-1\choose k}$ equals to the number of lattice paths in the $m\times n$ rectangle with $k$ external corners. Given a such path, let us continue it periodically to get an infinite path.
This construction maps exactly $k+1$ different paths (we have $k+1$ corners, as we have a corner at the starting point) into one. On the other hand, exactly one of them is totally above the diagonal - it corresponds to the set of corners with the lowest value of the linear function $my-nx$.
\end{proof}

\begin{lemma}
\label{Nara2}
$$\sum_{l=k}^{n-1}{l\choose k}N_{l}={(m+n-k-1)!\over m\cdot n\cdot (n-k-1)!(m-k-1)!k!}.$$
\end{lemma}
 
\begin{proof}
Remark that
$$N_{l}={k!(k+1)!(m-k-1)!(n-k-1)!\over l!(l+1)!(m-l-1)!(n-l-1)!}N_{k},$$
so
$${l\choose k}N_{l}={(k+1)!(m-k-1)!(n-k-1)!\over (l-k)!(l+1)!(m-l-1)!(n-l-1)!}N_{k}=$$ $${m-k-1\choose l-k}{n\choose n-l-1}{n\choose n-k-1}^{-1}N_{k}.$$
Now we have
$$\sum_{l=k}^{n-1}{m-k-1\choose l-k}{n\choose n-l-1}={m+n-k-1\choose n-k-1},$$
so
$$\sum_{l=k}^{n-1}{l\choose k}N_{l}={m+n-k-1\choose n-k-1}{n\choose n-k-1}^{-1}N_{k}={(m+n-k-1)!(k+1)!\over m!n!}N_{k}.$$
\end{proof}


\end{document}